\documentclass[a4j,10pt]{amsart}

\usepackage[mathscr]{eucal}
\usepackage{amscd}
\usepackage{amsfonts}
\usepackage{amsmath, amsthm, amssymb}
\usepackage{latexsym}
\usepackage[dvips]{graphics}
\usepackage{amssymb} 
\theoremstyle{plain} % イタリック体
\newtheorem{theorem}{\indent\sc Theorem}[section] % 見出しはスモールキャップ

\newtheorem{corollary}[theorem]{\indent\sc Corollary}
\newtheorem{proposition}[theorem]{\indent\sc Proposition}

\theoremstyle{definition} % ローマン体に変更
\newtheorem*{definition}{\indent\sc Definition}
\newtheorem*{conjecture}{\indent\sc Conjecture}

\newtheorem{remark}[theorem]{\indent\sc Remark}
\newtheorem{example}[theorem]{\indent\sc Example}

\def\C{{\mathbb{C}}}%   \C == \mathbb{C}
\def\R{{\mathbb{R}}}%   \R == \mathbb{R}
\def\Z{{\mathbb{Z}}}%   \Z == \mathbb{Z}
\def\D{{\mathbb{D}}}%  \D == \textbf{D}
\def\Pr{{\mathbb{P}}}

\title[Gauss map of minimal surfaces]{The Gauss map of pseudo-algebraic minimal surfaces}

\author[Y.~Kawakami, R.~Kobayashi, and R.~Miyaoka]{Y.~Kawakami, R.~Kobayashi, and R.~Miyaoka}

\thanks{The first author is partially supported by OCAMI (Osaka City university Advanced Mathematical Institute).}
\thanks{The second author is partially supported by Grant-in-Aid, 17204005, Ministry of Education and 
Science of Japan.
}
\thanks{The third author is partially supported by Grant-in-Aid, 16204007, Ministry of Education and 
Science of Japan.
}

\address{
Graduate school of Mathematics, 
Nagoya University, 
Nagoya, 464-8602
/Japan
}
\email{m02008w@math.nagoya-u.ac.jp}

\address{ % 第二著者
Graduate school of Mathematics, 
Nagoya University, 
Nagoya, 464-8602
/Japan
}
\email{ryoichi@math.nagoya-u.ac.jp}

\address{ % 第三著者
Graduate school of Mathematics, 
Kyusyu University, 
Fukuoka, 812-8581
/Japan
}
\email{r-miyaok@math.kyushu-u.ac.jp}
\begin{document}
\maketitle
\begin{abstract}
We refine Osserman's argument on the exceptional values of the Gauss map of 
algebraic minimal surfaces. This gives an effective estimate for the number of 
exceptional values and the totally ramified value number for a wider class of complete 
minimal surfaces that includes algebraic minimal surfaces. 
It also provides a new proof of Fujimoto's theorem for this class, 
which not only simplifies the proof but also reveals the
geometric meaning behind it. 
\end{abstract}
\section{Introduction}
The problem of finding the maximal number $D_{g}$ of the exceptional values of the Gauss map $g$ of 
a complete non-flat minimal surface $M$ in $\R^{3}$ was settled by Fujimoto \cite{F1}, \cite{F2} with the best 
possible upper bound being $4$. Indeed, for any number $r$, $0\leq r\leq 4$, we can construct complete 
minimal surfaces in $\R^{3}$ whose Gauss map omits exactly $r$ values. Moreover, Fujimoto proved that the 
totally ramified value number $\nu_{g}$, which gives more detailed information than $D_{g}$ does, satisfies 
$\nu_{g}\leq 4$, and this inequality is the best possible. Here, $b\in\Pr^{1}=\C\Pr^{1}$ is called 
a totally ramified value of $g\colon M\to \Pr^{1}$ if at all the inverse image points of $b$, $g$ branches. 
The exceptional values are regarded as totally ramified values, since it is natural to consider the multiplicity 
of an exceptional value to be infinite in the context of the Nevanlinna theory \cite{Ko}. The totally ramified 
value number $\nu_{g}$ is a weighted sum of the number of totally ramified values (see \S 3 for a precise definition). 
In particular, $D_{g}\leq \nu_{g}$ holds.

On the other hand, Osserman \cite{O1} proved that the Gauss map of a non-flat algebraic minimal surface 
omits at most $3$ values. By an algebraic minimal surface, we mean a complete minimal surface with 
finite total curvature. There are no known examples, however, of algebraic minimal surfaces whose Gauss map 
omits $3$ values, while there are many examples, of almost all topological types with the Gauss maps omitting 
$2$ values \cite{MS}. Thus an established conjecture is that the sharp upper bound of $D_{g}$ is $2$. 
Moreover, as in the case of Fujimoto's theorem, we have an implicit conjecture that the same is true for the 
totally ramified value number $\nu_{g}$. 

Surprisingly, the first author found algebraic minimal surfaces with totally ramified value number $\nu_{g}= 2.5$, 
i.e., strictly larger than $2$ \cite{Ka}. This overthrew the above implied conjecture, and a qualified problem 
to consider is then: 
\begin{center}
{\em Does there exist some $\kappa$, $2.5\leq \kappa< 3$ which is an upper bound $\nu_{g}\leq \kappa$ for $\nu_{g}$?} 
\end{center}

The totally ramified value number as well as the defect is well investigated in the Nevanlinna theory of transcendental 
meromorphic functions on $\C$. But for the Gauss map of {\em algebraic} minimal surfaces, this number has not 
been considered, at the best of the authors' knowledge. Since this is an indispensable number when we discuss 
the problem of exceptional values, we study it here by refining Osserman's algebraic argument, and the results 
turned out to be much more effective than we had expected (Theorem 3.3). They are effective in the sense that 
the upper bound we obtain is described in terms of the degree of the Gauss map and the topological data of $M$. 
Moreover, in some sense, this is the best possible result, and gives a proper extension of Osserman's results. 
In particular, we prove that the totally ramified value number of an algebraic minimal surface is strictly less than $4$. 

Another advantage of the approach here is that we can develop our arguments on a wider class of complete minimal surfaces, i.e., 
on all complete minimal surfaces whose Weierstrass data descends to meromorphic data on a compact Riemann surface. 
We refer to such minimal surfaces as {\em pseudo-algebraic minimal surfaces} (\S 3). As we do not assume the period condition, 
the surfaces may have infinite total curvature. For this class of minimal surfaces, we obtain a new proof of Fujimoto's 
theorem, which reveals the reason why $\nu_{g}\leq 4$ holds. We also give a kind of unicity theorem (Theorem 5.1), 
which asks the least number of values at which if two Gauss maps $g_{1}$ and $g_{2}$ have the same inverse image then $g_{1}=g_{2}$.

Moreover, our argument suggests how to estimate the characteristic function $T_{g}(r)$ of the Gauss map lifted 
to the universal covering surface, which plays an essential role in the Nevanlinna theory on the unit disk (\S 6).  
We believe that our results give an important link to the Nevanlinna theory for future research. 

\section{Preliminaries}
A minimal surface $x\colon M\to \R^{3}$ is originally considered as a surface spanning a given frame 
with least area. The Euler-Lagrange equation turns out to be 
\begin{equation}\label{harmonic}
\triangle x=0\,,
\end{equation}
i.e., each coordinate function of a minimal surface is harmonic.
Thus there exist no compact minimal surfaces
without boundary.
With respect to a complex parameter $z=u+iv$ of the surface,
(\ref{harmonic}) is given by $\bar\partial\partial x=0$ 
where $\partial=\dfrac{\partial}{\partial u}-i\dfrac{\partial}{\partial v}
$ following Osserman \cite{O2}. Thus if we put 
$$
\partial x=(\phi_1,\phi_2,\phi_3)\,,
$$
$\phi_j$'s are holomorphic differentials on $M$.
These satisfy 
\begin{quote}
(C) $\sum \phi_{i}^{2}=0$ : conformality condition \\
(R) $\sum |\phi_{i}|^{2}>0$ : regularity condtion \\
(P) For any cycle $\gamma \in H_{1}(M,\Z)$, $\Re\int_{\gamma}\phi_{i}=0$ : period condtion \\
\end{quote}
We recover $x$  by
the real Abel-Jacobi map (called the Weierstrass-Enneper representation 
formula) 
\begin{equation}\label{EWrep.}
x(z)=\Re\int_{z_0}^z(\phi_1,\phi_2,\phi_3)
\end{equation}
up to translation. 
From here on we restrict our attention to non-flat minimal surfaces. 
If we put
\begin{equation}\label{W-data}\left\{
\begin{array}{ll}
hdz=\phi_1-i\phi_2\\
g=\dfrac{\phi_3}{\phi_1-i\phi_2}
\end{array}\right.
\end{equation}
then $hdz$ is a holomorphic differential and $g$ is a
meromorphic function on $M$.
Geometrically, it is well-known that $g$ is the stereographically
projected Gauss map of $M$.
We call $(hdz,g)$ the Weierstrass data. 
This is related to $\phi_j$'s\,  in a one to one way by
\begin{equation}\label{phis}\left\{
\begin{array}{ll}
\phi_1=\dfrac{h}{2}(1-g^2)dz\\
\phi_2=\dfrac{ih}{2}(1+g^2)dz\\
\phi_3=hgdz\,.
\end{array}\right.
\end{equation}
If we are given a holomorphic differential $hdz$ and a
meromorphic function $g$ on $M$, we get $\phi_j$'s\, by
this formula. They satisfy (C) automatically, and 
the condition (R) is interpreted as the poles of $g$ of order $k$
coincides exactly with the zeros of $hdz$ of order $2k$, because
the induced metric on $M$ is given by 
$$
ds^2=\frac{|h|^2(1+|g|^2)^2}{4}|dz|^2.
$$
A minimal surface is complete if all divergent paths have 
infinite length with respect to this metric. 
In general, for a given meromorphic function $g$ on $M$, 
it is not so hard to find a holomorphic differential
$hdz$ satisfying (R).
But the period condition (P) always causes trouble.
When (P) is not satisfied, we anyway obtain a minimal surface on the 
universal covering surface of $M$.

Here we notice that the triple of holomorphic differentials
$$
e^{i\theta}(\phi_1,\phi_2,\phi_3),\quad \theta\in \R
$$
also satisfies (C) and (R).
The corresponding  Weierstrass data is given by
\begin{equation}\label{associated}
\left\{
\begin{array}{ll}
g^\theta(z)=g(z)\\
h^\theta dz=e^{i\theta}hdz\,.
\end{array}
\right.
\end{equation}
As (P) is scarcely satisfied by these data,
we get an  $S^1$ parameter family of minimal surfaces defined 
on the universal covering surface by (\ref{EWrep.}), which is called 
the associated family.
Note that all surfaces in this family have the same 
Gauss map.

\begin{example}
Catenoid and Helicoid are well-known as surfaces 
belonging to the same associated family. 
The Gauss map of this family omits two values.
\end{example}

Now the Gauss curvature $K$ of $M$ is given by
\begin{equation}
K=-\Bigl(\frac{4|g'|}{|h|(1+|g|^2)^2}\Bigr)^2
\label{curvature}
\end{equation}
and the total curvature by
\begin{equation}\label{totalcurv.}
\tau(M)=\int_MKdA=-\int_M \Bigl(\frac{2|g'|}{(1+|g|^2)}\Bigr)^2du\wedge dv
\end{equation}
where $dA$ is the surface element of $M$.
Note that $|\tau(M)|$ is the area of $M$ with respect to
the (singular) metric induced from the Fubini-Study metric of 
$\Pr^1$ by $g$. 
When the total curvature of a complete minimal surface is finite,
the surface is called {\em an algebraic minimal surface}.

\begin{theorem}[Huber, Osserman]\label{Huber}
An algebraic minimal surface
$x:M\to \R^3$ satisfies: 
\begin{enumerate}
\item[\rm(i)] $M$ is conformally equivalent to 
$\overline{M}\setminus\{p_1,\dots,p_k\}$ where $\overline{M}$ 
is a compact Riemann surface, and $p_1,\dots,p_k\in\overline{M}$ \: \cite{H}.
\item[\rm(ii)]  The Weierstrass data $(hdz,g)$ is extended 
meromorphically to $\overline{M}$ \:\cite{O1}.
\end{enumerate}
\end{theorem}

We denote the number of exceptional values of $g$ by $D_g$.
Other than Catenoid, there are many examples of algebraic minimal 
surfaces with $D_g=2$,  which include those of hyperbolic type. 

\begin{theorem}[Miyaoka-Sato \cite{MS}]\label{MS1}
There exist algebraic minimal surfaces with $D_g=2$, for
\begin{enumerate}
\item[\rm{(i)}] $G=0, k\geq 2$
\item[\rm{(ii)}]  $G=1, k\geq 3$
\item[\rm{(iii)}]  $G\geq 2, k\geq 4$
\end{enumerate}
where $G$ $($resp. $k$$)$ is the genus $($resp. the number of
punctures$)$ of the Riemann surface on which the surfaces are defined. 
\end{theorem}

When $G=0$ and $k=2$, all such minimal surfaces are classified.
Examples for $G=0$ and $k=3$  
given below [MS, Proposition 3.1] are important for later argument:
 let $M=\Pr^1\setminus\{\pm i,\infty\}$,
and define a Weierstrass data by
\begin{equation}\label{MS2}
\left\{
\begin{array}{ll}
g(z)=\sigma\dfrac{z^2+1+a(t-1)}{z^2+t}\\
hdz=\dfrac{(z^2+t)^2}{(z^2+1)^2}dz,\quad (a-1)(t-1)\ne 0\\
\sigma^2=\dfrac{t+3}{a\{(t-1)a+4\}}\,.
\end{array}
\right.
\end{equation}
For any $a,t$ satisfying $\sigma^2<0$, we obtain an algebraic
minimal surface of which Gauss map omits two values 
$\sigma,\sigma a$.

Applying the covering method to this surface 
(see Remark 3.6), we obtain examples 
of (ii) and (iii).
But as these examples have all the same image in $\R^3$, we further 
 constructed mutually non-congruent examples for $G=1$ and $k=4$, by generalizing Costa's surface 
[MS,Theorem 3].  For details see Remark 4.2.

\section{Pseudo-algebraic minimal surfaces and the main results}
\begin{definition}
We call a complete minimal surface in $\R^3$ {\em pseudo-algebraic}, 
if the following conditions are satisfied:
\begin{enumerate}
\item[(i)] The Weierstrass
data $(hdz,g)$ is defined on a Riemann surface $M=\overline{M}\setminus\{p_1,\dots,p_k\}$, $p_j\in \overline{M}$, 
where $\overline{M}$ is a compact Riemann surface.
\item[(ii)] $(hdz,g)$ can be extended meromorphically to $\overline{M}$.
\end{enumerate}
We call $M$ {\em the basic domain} of the pseudo-algebraic minimal surface. 
\end{definition}

Since we do not assume the period condition on $M$, 
a pseudo-algebraic minimal surface is defined on some covering surface of $M$, 
in the worst case, on the universal covering.
Note that Gackst\"atter called such surfaces {\em abelian minimal surfaces} \cite{G}.

Algebraic minimal surfaces and their associated surfaces are
certainly pseudo-algebraic. Another important example is Voss' surface.
The Weierstrass data of this surface is defined on $M=\C\setminus\{a_1,a_2,a_3\}$ for
distinct $a_1,a_2,a_3\in \C$, by
\begin{equation}\label{Voss}
\left\{
\begin{array}{ll}
g(z)=z\\
hdz=\dfrac{dz}{\Pi_j(z-a_j)}\,.
\end{array}
\right.
\end{equation}
As this data does not satisfy the period condition,
we get a minimal surface $x:\D\to \R^3$ on the universal
covering disk of $M$.
In particular, it has infinite total curvature.
We can see that the surface is complete
 and the Gauss map omits four values $a_1,a_2,a_3,\infty$.
Starting from $M=\C\setminus\{a_1,a_2\}$, we get similarly 
a complete minimal surface $x:\D\to \R^3$, 
of which Gauss map omits three values $a_1,a_2,\infty$.
The completeness restricts the number of points $a_j$'s 
to be less than four. 
Note that in  both cases, no elements of the associated 
family satisfy (P), 
hence  have infinite total curvature.

\begin{remark}\label{Lopez}
There exist complete minimal surfaces with $D_{g}=4$ which are not pseudo-algebraic
 (see \cite{L}).
\end{remark}
Now we define the totally ramified value number $\nu_{g}$ of $g$.
\begin{definition}
We call $b\in \Pr^1$ a totally ramified value of $g$ when at any 
inverse image of $b$, $g$ branches. 
We regard exceptional values also as totally ramified values.
Let $\{a_1,\dots,a_{r_0},b_1,\dots,b_{l_0}\}\subset \Pr^1$ be the set of 
 totally ramified values of $g$, where $a_j$'s are exceptional 
values. For each $a_j$, put $\nu_j=\infty$, and for each $b_j$, 
define $\nu_j$ to be the minimum of the multiplicity of $g$ at 
points $g^{-1}(b_j)$. Then  we have $\nu_j\geq 2$.
We call 
\[
\nu_{g}=\sum_{a_j,b_j}\biggl(1-\dfrac1{\nu_j}\biggl)=r_0+\sum_{j=1}^{l_0}\biggl(1-\dfrac1{\nu_j}\biggl)
\]
{\em the totally ramified value number of $g$}.
\end{definition}
To explain the natural meaning of this number, 
we need the second main theorem in the Nevanlinna theory, 
 which we have no space to mention here (See \cite{Ko}). 
Note that though $\nu_{g}$ is a rational number, 
the upper bound is given 
by the integer 4 by Fujimoto's theorem.

\begin{theorem}[Kawakami \cite{Ka}]\label{2.5}
The Gauss map of the algebraic minimal surfaces 
given in $(\ref{MS2})$ 
has totally ramified value number 2.5. 
\end{theorem}

In fact, it has two 
exceptional values, and another totally ramified value at $z=0$ where 
$g'(z)=0$.
This theorem is a breakthrough to propose 

\begin{conjecture}
For {\em algebraic} minimal surfaces,  
there exists $2.5\leq \kappa<3$ which satisfies $\nu_{g}\le \kappa$.
\end{conjecture}

Note that the conjecture is {\em not} for pseudo-algebraic minimal
 surfaces, since Voss' surfaces satisfy $\nu_{g}=3$ and 4.
But we 
develop an algebraic argument on the pseudo-algebraic minimal surfaces 
in a unified way,  
and then specialize the results in the algebraic case.  
The following is the main result of this paper.

\begin{theorem} \label{main1}
Consider a non-flat pseudo-algebraic minimal surface with the basic domain $M=\overline{M}\setminus\{p_1,\dots,p_k\}$.
Let $G$ be 
the genus of $\overline{M}$, and let $d$ be the degree of $g$
considered as a map on $\overline{M}$.
Then we have
\begin{equation}\label{estimate1}
D_g\leq 2+\dfrac2{R}\,,\quad
\dfrac{1}{R}=\dfrac{G-1+k/2}{d}\leq 1\,.
\end{equation}
More precisely, if the number of (not necessarily totally) 
ramified values other than the exceptional values of $g$ is $l$, 
we have 
\begin{equation}\label{estimate2}
D_g\leq 2+\dfrac2{R}-\dfrac{l}{d}\,.
\end{equation}
On the other hand, the totally ramified value number of $g$  
satisfies 
\begin{equation}\label{estimate3}
 \nu_{g}\leq 2+\dfrac2{R}\,. 
\end{equation}
In particular, we have
\begin{equation}\label{estimate4}
D_g\leq \nu_{g}\leq 4\,,
\end{equation}
and for algebraic minimal surfaces, the second inequality is a 
strict inequality. 
(\ref{estimate2}) and (\ref{estimate3}) are best possible in both 
algebraic and non-algebraic cases. 
\end{theorem}

The geometric meaning of the ratio $R$ is given in \S 6. 
This theorem implies the following known facts:
\begin{corollary} [Osserman, Fang, Gackst\"atter] \label{main2}
For algebraic minimal surfaces, we have:
\begin{enumerate}
%\item $D_g\leq 3$, i.e., Osserman's theorem.
\item[\rm(i)] When $G=0$, $D_g\le 2$ holds.
\item[\rm(ii)] When $G=1$ and $M$ has a non-embedded end, $D_g\le 2$ holds.
If $G=1$ and $D_g=3$ occur, $d=k$ follows and $g$ does not branch 
in $M$, so is a non-branched covering of $\Pr^1\setminus 
\{3\text{ points}\}$.
\end{enumerate}
\end{corollary}
\begin{remark}
Fang [Fa, Theorem 3.1] shows that algebraic 
minimal surfaces with $d\leq 4$ satisfy $D_g\leq 2$  (see \cite{WX}
for $d\leq 3$).
\end{remark}
\begin{remark} 
There exists a way of construction of algebraic minimal surfaces 
by a covering method of Klotz-Sario \cite{BC}.
Indeed, if $x:M\to \R^3$ is an algebraic minimal surface,
and if $\pi:\hat M\to M$ is a non-branched covering surface 
of $M=\overline{M}\setminus\{p_1,\dots,p_k\}$,
then we obtain a new algebraic minimal surface by 
$\hat x=x\circ \pi:\hat M\to \R^3$. 
This surface has the same image as the original one, but
the domain $\hat M$ has different topological type.
Nevertheless, we can see that {\em the ratio $R$ is invariant}
under this construction, via a little algebraic argument.
Certainly, $D_g$ and $\nu_{g}$ are also invariant under covering
construction.
\end{remark}
\section{Proof}
The proof of Theorem \ref{main1} is given by a refinement of the proof 
of Osserman's theorem in \cite{O1}. In order to simplify the argument, 
we may assume without loss of generality 
that $g$ is neither zero nor pole at  $p_j$, and moreover, 
zeros and poles of $g$ are simple.
By completeness, $hdz$ has poles of order $\mu_j\ge 1$ at $p_j$.
The period condition implies $\mu_j\geq 2$, 
but here we do not assume this.
Let $\alpha_s$ be (simple) zeros of $g$, $\beta_t$ (simple) poles of $g$.
The following table shows the relation between zeros and poles 
of $g$, $hdz$ and $ghdz$. 
The upper index means the order. 
\begin{center}
\begin{tabular}{|c|c|c|c|}\hline
$z$ & $\alpha_s$  & $\beta_t$ & $p_j$  \\\hline
$g$& $0^1$ & $\infty^1$ &   \\\hline
$hdz$ &   & $0^2$ & $\infty^{\mu_j}$  \\\hline
$ghdz$ & $0^1$ & $0^1$ & $\infty^{\mu_j}$ \\\hline
\end{tabular}
\end{center}
Applying the Riemann-Roch formula to the meromorphic
differential $hdz$ or $ghdz$ on $\overline{M}$, we obtain
\[
2d-\sum_{j=1}^k\mu_j=2G-2\,.
\]
Note that this equality depends on the above setting of
zeros and poles of $g$, though $d$ is an invariant.
Thus we get
\begin{equation}\label{cowen}
d=G-1+\dfrac12\sum_{j=1}^k\mu_j\ge G-1+\dfrac{k}2\,,
\end{equation}
and
\begin{equation}\label{ration}
R^{-1}\leq 1\,.
\end{equation}
When $M$ is an algebraic minimal surface or its associated surface,
we have $\mu_j\ge 2$ and so $R^{-1}<1$.

Now, we prove (\ref{estimate2}) (and (\ref{estimate3})).
Assume $g$ omits $r_0=D_g$ values, and let $n_0$ be the sum of the 
branching orders of $g$ at these exceptional values.
Moreover, let $n_b$ be the sum of branching orders at the inverse
images of  
non-exceptional (not necessarily totally) ramified values $b_1,\dots,b_l$ of $g$.  
We see
\begin{equation}\label{puncture}
k\ge dr_0-n_0\,,\quad 
n_b\ge l\,.
\end{equation}
Let $n_{g}$ be the total branching order of $g$. 
Then applying Riemann-Hurwitz's theorem to the 
meromorphic function $g$ on $\overline{M}$,
we obtain
\begin{equation}\label{Riemann-Hurwitz}
n_{g}=2(d+G-1)=n_0+n_b\ge dr_0-k+l\,.
\end{equation}
If we denote
\[
\nu_i=\text{min}_{g^{-1}(b_i)}\{\text{multiplicity of }g(z)=b_i\}\,,
\]
we have $1\le \nu_i\le d$. 
Now the number of exceptional values satisfies
\begin{equation}
D_g= r_0\leq \dfrac{n_{g}+k-l}{d}=2+\dfrac2{R}-\dfrac{l}{d}
\end{equation}
where we have used (\ref{Riemann-Hurwitz}), hence (\ref{ration}) implies
\[
D_g\leq 2+\dfrac2{R}\leq 4\,.
\]
In particular for algebraic minimal surfaces and its associated 
family, we have $R>1$ so that
\[
D_g\le 3\,,
\]
which is nothing but Osserman's theorem.

Next, we show (\ref{estimate4}).
Let $b_1,\dots,b_{l_0}$ be the {\em totally} ramified values which are not exceptional values. 
Let $n_r$ be the sum of branching orders 
at $b_1,\dots,b_{l_0}$. For each $b_i$, the number of points
in the inverse image $g^{-1}(b_i)$ 
is less than or equal to $d/\nu_i$, since $\nu_i$ is the 
minimum of the multiplicity at all $g^{-1}(b_i)$.
Thus we obtain
\begin{equation}\label{TRVN}
dl_0-n_r\leq \sum_{i=1}^{l_0}\dfrac{d}{\nu_i}\,.
\end{equation}
This implies
\[
l_0-\sum_{i=1}^{l_0}\dfrac{1}{\nu_i}\leq \dfrac{n_r}{d}\,,
\]
hence using the first inequality in (\ref{puncture}) and $n_r\leq n_b$, we get
\[
\nu_{g}=r_0+\sum_{i=1}^{l_0}\biggl(1-\dfrac{1}{\nu_i}\biggl)\leq 
\dfrac{k+n_0}{d}+\dfrac{n_r}{d}\leq \dfrac{n_{g}+k}{d}=2+\dfrac2{R}\,.
\]
The sharpness of (\ref{estimate2}) and (\ref{estimate3}) follows from:
\begin{enumerate}
\item[(1)] When $d=2$ we have
\[
D_g\le 2+\frac2{R}-\frac{l}{2}, \quad \nu_{g}\leq 2+\dfrac2{R}\,.
\]
The surface given by (\ref{MS2}) attains both equalities, since $R=4$, $l=1$ and
$D_g=2$, $\nu_{g}=2.5$. Thus (\ref{estimate2}) and (\ref{estimate3}) are
sharp.
\item[(2)] 
Voss' surface satisfies $d=1$ and $G=0$.
Thus when $k=3$, we get $R=2$, $l=0$ hence 
$D_g=3=2+2/2$.
When $k=4$, we have $R=1$, $l=0$ and 
$D_g=4=2+2/1$. These show that (\ref{estimate2}) and (\ref{estimate3}) are
sharp in non-algebraic pseudo-algebraic case, too.
\end{enumerate}

\bigskip%\noindent
Corollary 3.4 is obtained as follows:
 It is easy to see that $r_0=3$ implies $R\leq 2$, hence 
$G-1+\dfrac12\sum_{j=1}^k\mu_j\leq 2(G-1)+k$.
As we have $\mu_j\geq 2$ in the algebraic case, it follows
\begin{equation}\label{embedded1}
k\leq \dfrac12\sum_{j=1}^k\mu_j\leq G-1+k\,.
\end{equation}
Thus  we obtain (i). 
When $G=1$, (\ref{embedded1}) implies $\mu_{j}=2$ for all $j$, which means that all 
the ends are embedded (\cite{JM}). Therefore, if $M$ has a non-embedded end, then $r_{0}\leq 2$. 
When $r_{0}=3$, from (\ref{cowen}), $d=k$ holds, and hence $R=2$. Therefore, from (\ref{estimate2}), 
$l=0$ holds, which means that $g$ does not branch in $M$.
\qed

\begin{remark} 
The inequality (\ref{estimate2})  gives
more informations than (\ref{estimate1}). 
In particular, (\ref{estimate2}) implies that the more branch points 
$g$ has in $M$,  the less is the number of exceptional 
values. 
\end{remark}
\begin{remark} The inequality (\ref{estimate1}) is also best possible 
for algebraic minimal surfaces in the following sence. 
In [MS, Theorem 3], we constructed two infinite 
series of mutually distinct algebraic minimal surfaces 
of the fixed topological type $G=1$ and $k=4$, 
whose Gauss map omits 2 values.  
(There are some errors in signatures in [MS, Lemma 4.1], but 
no effects on the result.)
They are given as follows. Let $\overline{M}$ be the square torus  
on which the Weierstrass $\wp$ 
function satisfies $(\wp')^2=4\wp({\wp}^2-a^2)$. 
Let $M$ be given by removing 
4 points satisfying $\wp=0,\pm a,\infty$ from $\overline{M}$.
Define the Weierstrass data by

\medskip
\noindent[Case 1] $g=\dfrac\sigma{{\wp}^j{\wp}'},\quad 
hdz=\dfrac{{\wp}d{\wp}}{{\wp}'},\quad j=1,2,3\dots,$

\noindent[Case 2] $g=\dfrac\sigma{{\wp}^j{\wp}'},\quad 
hdz=\dfrac{{\wp}^{j+1}d{\wp}}{{\wp}'},\quad j=2,4,6\dots,$

\medskip\noindent
Then choosing a suitable $\sigma$, we obtain algebraic minimal 
surfaces with $g$ omitting 2 values 0 and $\infty$.
Since the degree of $g$ is $d=2j+3$ in both cases and  
$R=d/2=(2j+3)/2$,  
$2+2/R$ tends to $2$ $(=D_g)$ as close as we like. 
(Costa's surface is given by $j=0$, in which case
$(G,k,d)=(1,3,3)$, and $g$ omits just one value 0.)
\end{remark}

\section{Unicity theorem and some other results}
We give two applications of Theorem \ref{main1}. 
The first one is an extension of Fujimoto's unicity theorem for algebraic minimal surfaces \cite{F3} 
to the pseudo-algebraic case. 
By using Fujimoto's argument and Theorem \ref{main1}, 
we obtain the following result:
\begin{theorem}\label{main3}
Consider two non-flat pseudo-algebraic minimal surfaces $M_1, M_2$ with the 
same basic domain $M=\overline{M}\setminus\{p_1,\dots,p_k\}$.
Let $G$ be the genus of $\overline{M}$, and let 
$g_1, g_2$ be the Gauss map of $M_1$ and $M_2$ respectively.
Assume that $g_1$ and $g_2$ have the same degree $d$ 
when considered as a map on $\overline{M}$, but assume $g_1\ne g_2$ 
as a map $M\to \Pr^1$. Let $c_1,\dots, c_q\in \Pr^1$ be
distinct points such that $g_1^{-1}(c_j)\cap M=g_2^{-1}(c_j)\cap M$ for 
$1\leq j\leq q$. Then 
\begin{equation}\label{unicity1}
q\leq 4+\frac2{R},\quad \dfrac{1}{R}=\dfrac{G-1+k/2}{d}\leq 1
\end{equation}
follows. In particular, $q\leq 6$, and for algebraic 
minimal surfaces we have $q\leq 5$.
\end{theorem}
\begin{proof}
Put
\[
\delta_j=\sharp(g_1^{-1}(c_j)\cap M)=\sharp (g_2^{-1}(c_j)\cap M)\,,
\]
where $\sharp$ denotes the number of points.
Then we have
\begin{equation}\label{unicity2}
qd\leq k+\sum_{j=1}^q\delta_j+n_{g}\,,
\end{equation}
using the same notation as in \S 4.
Consider a meromorphic function $\varphi=\dfrac1{g_1-g_2}$ on $M$.
Then at each point of $g_1^{-1}(c_j)\cap M$, $\varphi$ has 
a pole, while the total number of the poles of $\varphi$ 
is at most $2d$, hence we get
\begin{equation}\label{unicity3}
\sum_{j=1}^q\delta_j\leq 2d\,.
\end{equation}
Then from (\ref{unicity2}) and (\ref{unicity3}), we obtain
\[
qd\leq k+2d+n_{g}\,,
\]
and
\[
q\leq \dfrac{2d+n_{g}+k}{d}=4+\dfrac2{R}
\]
follows immediately.
\end{proof}
\begin{remark}
Fujimoto \cite{F3} gives an example of two 
pseudo-algebraic minimal surfaces with $q=6$, of which Gauss maps
do not coincide. For algebraic case, whether $q=5$ is best possible or
not is another interesting open problem.
\end{remark}
%\medskip
%\Remark The degree condition in the theorem can be omitted. 
%In fact, if $d_i$ is the degree of $g_i$, then let $d$ be 
%the LCM of $d_1$ and $d_2$,
%and take the covering $\hat M_i$ of each surface so that the degree 
%of the Gauss map becomes $d$. Then apply above argument to
%$\hat M_1$ and $\hat M_2$. 

The second application of Theorem \ref{main1} is a proof of Gackst\"atter's result \cite{G} :
\begin{proposition}[Gackst\"atter \cite{G}]\label{Gack} 
If the Gauss map of an algebraic minimal surface with $G=1$ 
%and its associated surfaces has a Gauss map with $D_g=3$, 
omits 3 values $a_1,a_2,a_3\in \Pr^1$, 
then all branch points of $g$ are located at the end points, and
$g$ is  a non-branched covering map 
 of $\Pr^1\setminus\{a_1,a_2,a_3\}$. 
\end{proposition}
\begin{proof}
This follows immediately from Corollary 3.4 (ii).  
\end{proof}

Thus the Gauss map descends to 
$\Pr^1\setminus \{3\text{ points}\}$, but the minimal surface 
is not obtained from a covering of a minimal surface 
defined on $\Pr^1\setminus \{3\text{ points}\}$,
otherwise, by (ii) of Corollary 3.4, $D_g\leq 2$. 
This implies that $hdz$ can not descends to 
$\Pr^1\setminus \{3\text{ points}\}$.

We will use Proposition \ref{Gack} in order to give a lower bound for the 
characteristic function $T_g(r)$ in Proposition 6.3. 

The following is obvious:
\begin{proposition}
If the Gauss map $g$ of a pseudo-algebraic minimal surface 
omits exactly $r$ values $a_1,\dots,a_r\in \Pr^1$ for $r=3,4$, and has 
no branch points in the basic domain $M$, then $g$ is 
a non-branched covering of $\Pr^1\setminus \{a_1,\dots,a_r\}$.
\end{proposition}

In this situation, the universal covering 
surface of $M$ and of $\Pr^1\setminus\{a_1,\dots,a_r\}$ are disks,
which we denote by $\D$ and $\Omega$, respectively.
When $g$ has no branch points in $M$, 
 the lifted map
$g:\D\to \Omega$ is a non-branched holomorphic map, 
i.e., a hyperbolic isometry. 
Since the degree of $g$ restricted to $\overline{M}$ is $d$, 
the fundamental domain of $M$ is given by $\cup_{i=1}^d U_i\subset \D$,  where each $U_i$ is diffeomorphic to 
$\Pr^1\setminus \{a_1,\dots,a_r\}$. 
\begin{example}
Voss' surfaces are examples for $d=1$.
\end{example}
\section{Toward the Nevanlinna theory}
Unfortunately, the above argument does not prove the conjecture in \S 3. 
To go further, we state some links to the Nevanlinna theory 
\cite{Ko}. 
 
We consider the case where the universal covering 
surface of $M$ is a unit disk $\D$.
In order to adjust to the Nevanlinna theory, we use the 
 hyperbolic metric $\omega_h$ with curvature 
$-4\pi$ on $\D$, and 
the Fubini-Study metric $\omega_{FS}$ with curvature $4\pi$ on $\Pr^1$ 
(hence $\Pr^1$ has area 1).
Then by Gauss-Bonnet's theorem for a complete punctured Riemann
surface with hyperbolic metric, we have
\begin{equation}\label{GaussBonnet}
2\pi\chi(M)=\int_M K_h \omega_h=-4\pi \int_M \omega_h
=- 4\pi A_{hyp}(M)\,,
\end{equation}
where $A_{hyp}(M)$ is 
the hyperbolic area of $M$, hence for the fundamental domain $F$ of 
$M$, we get
\begin{equation}\label{hyperbolicarea}
A_{hyp}(F)=G-1+\dfrac{k}{2}\,.
\end{equation}
\begin{remark}
The Gauss-Bonnet theorem (\ref{GaussBonnet}) for $(M,\omega_h)$ is often used without proof, 
so here we give a brief proof. 
Let $D_{\varepsilon_j}$
be the disk with radius $\varepsilon_j$ around $p_j$, $j=1,2,\dots, k$.
We denote $M_{\varepsilon}=\overline{M}\setminus \cup_jD_{\varepsilon_j}$,
and by $\varepsilon\to 0$, we mean all $\varepsilon_j\to 0$.
Consider any metric $\sigma$ on $\overline{M}$ which is flat in all $D_{\varepsilon_j}$. Denoting locally (as K\"ahler forms) 
$\sigma=\dfrac{i}2\tilde \sigma dz\wedge d\bar z$  
and $\omega_h=\dfrac{i}2\tilde \omega_hdz\wedge d\bar z$, 
we have by Stokes' theorem 
\[
-\sum_j\int_{\partial D_{\varepsilon_j}}d^c\log (\sigma/\omega_h)
=\int_{M_{\varepsilon}}dd^c\log (\sigma/\omega_h)
=\int_{M_{\varepsilon}}dd^c\log \tilde\sigma
-\int_{M_{\varepsilon}}dd^c\log \tilde\omega_h\,,
\]
where $d=\partial+\bar\partial$, $d^c=(\partial-\bar\partial)/(4\pi i)$, 
(here $\partial$ is the half of Osserman's one). 
Because $dd^c\log\tilde\omega=-\dfrac{K_\omega}{2\pi}dA_{\omega}$ 
holds where $K_\omega$ and $dA_\omega$ are the Gauss curvature 
and the area form of $\omega$, respectively, 
taking the limit $\varepsilon\to 0$ and applying 
the Gauss-Bonnet's theorem to $(\overline{M},\sigma)$, we obtain  
\[
\lim_{\varepsilon\to 0}
\int_{M_{\varepsilon}}(dd^c\log \tilde\sigma-dd^c\log \tilde\omega_h)
=-\chi(\overline{M})-2A_{hyp}(M)\,.
\]
Next, take a local coordinate on each $D_{\varepsilon_j}$ 
so that $z=0$ corresponds to $p_j$. Then we can express 
$\sigma=\dfrac{i}2dz\wedge d\bar z$ 
and
 $\omega_h=\dfrac{i}{2\pi}\dfrac{dz\wedge d\bar z}{|z|^2(\log|z|^{-2})^2}$ on $D_{\varepsilon_j}$.
Noting that $d^c=\dfrac1{4\pi}\bigl(-\dfrac1{r}\dfrac{\partial}{\partial\theta}dr+r\dfrac{\partial}{\partial r}d\theta\bigr)$, we obtain 
\[
\lim_{\varepsilon\to 0}\sum_j\int_{\partial D_{\varepsilon_j}} d^c\log (\sigma/\omega_h)=k\,,
\]
which implies (\ref{GaussBonnet}) and (\ref{hyperbolicarea}).
\qed
\end{remark}
Next, let $d$ be the degree of $g$, 
then the area $A_{FS}(F)$ of $F$
with respect to the induced metric $g^*\omega_{FS}$ is  $d$.
Thus we obtain
\begin{equation}\label{ration2}
A_{FS}(F)=\dfrac{d}{G-1+k/2}A_{hyp}(F)=RA_{hyp}(F)\,.
\end{equation}
We now know the meaning of the ratio $R$;  the ratio of 
the area of the fundamental domain with respect to the 
induced Fubini-Study metric
to the one with respect to the hyperbolic metric on $\D$.
\begin{remark}\label{nothyper}
Even when the conformal type of $M$ is not  
hyperbolic, the ratio $R$ is meaningful in 
Theorem \ref{main1}.
\end{remark}
Now, remember Shimizu-Ahlfors' theorem on the characteristic
function $T_g(r)$ of $g$, in terms of
\[
T_g(r)=\int_0^r \dfrac{dt}{t}\int_{\C(t)}g^*\omega_{FS}\,.
\]
Here $\C(t)$ is the subdisk of $\D$ 
with radius $0<t<1$. 
In order to develop the Nevanlinna theory on meromorphic functions 
on the unit disk, we need the growth order of  $T_g(r)$ 
compared with
\[
\int_0^r\dfrac{dt}{t}\int_{\C(t)}\omega_{h}\approx
\dfrac12\log\dfrac1{1-r}\,,
\]
where $r$ is sufficiently close to 1
(strictly, the left hand side is
$\dfrac12\log\dfrac1{1-r^2}$). 
We always use this approximation formula
in the following, because in the Nevanlinna theory,
a bounded quantity is ignored.

Actually, we want to know the best estimate of type
\begin{equation}\label{Simizu1}
T_g(r)\geq \dfrac{\eta}{2}\log\dfrac1{1-r}\,,
\end{equation}
to get the lemma on logarithmic derivatives on meromorphic functions on 
the unit disk  \cite{KKM}. 
In general, the area of $\C(t)$ is approximated by that of finite
union of fundamental domains $\cup F_j$. 
It seems that we get $\eta$ from (\ref{ration2}),
but we need some argument here, 
since the hyperbolic symmetry never fits the shape 
of the disk. 
We do not go into details of this argument, instead,
we give some examples in easier cases.
\newline

If we replace the Fubini-Study 
metric by a singular metric $\Psi$ on $\Pr^1$ with area 1, 
we have
\begin{equation}\label{Simizu2}
T_g(r)\geq \displaystyle\int_0^r\dfrac{dt}{t}\int_{\C(t)}g^*\Psi\,.
\end{equation}
This is shown rather easily by using Crofton's formula 
in the integral geometry \cite{Ko}.
When the image $g(M)$ is $\Pr^1\setminus \{r\text{ points}\}$,
where $r=3$ or $4$, the singular metric $\Psi$ on $\Pr^1$ 
induced by the hyperbolic metric on $\Omega$ normalized so that the 
area of $g(M)$ (counted without multiplicity) is $1$ fits the case.
Using this metric, we give a few computable examples. 
\begin{proposition}\label{Simizu3}
Consider a pseudo-algebraic minimal surface with 
the basic domain $M=\overline{M}\setminus\{p_1,\dots,p_k\}$, 
and assume that $g$ branches only at $p_j$'s.   
\begin{enumerate}
\item[\rm{(i)}] If $D_g=3$, we have 
\begin{equation}\label{case3}
T_g(r)\geq \log\dfrac1{1-r}\,.
\end{equation}
This is satisfied by Voss' surface with $k=3$, 
and an algebraic minimal surface with $G=1$ and $D_g=3$, if any.
\item[\rm{(ii)}] If $D_g=4$, we have 
\begin{equation}\label{case4}
T_g(r)\geq \dfrac12\log\dfrac1{1-r}\,.
\end{equation}
This is satisfied by Voss' surface with $k=4$. 
\end{enumerate}
\end{proposition}
\begin{remark}
For some reasons, we conjecture that 
(\ref{case3}) and (\ref{case4}) are
equalities, and  $\eta=2$ in (i) and $\eta=1$ in (ii) hold, 
which implies $\eta=R$ in these cases. 
\end{remark}
\begin{proof}
Let $\D$ be the universal covering disk of $M$, and 
$\Omega$ that of $\Pr^1\setminus\{a_1,\dots,a_{r_0}\}$, 
where $a_1,\dots,a_{r_0}$ are the exceptional values of $g$.
Let $\omega_{h}$ and $\omega_{\Omega}$ be the 
hyperbolic metric with curvature $-4\pi$.
Denote by $g:\D\to \Omega$ the lifted 
Gauss map. Since this is not branched, $g$
is a hyperbolic isometry.
To obtain the characteristic function $T_g(r)$, 
normalize the metric $\omega_{\Omega}$ so that 
the fundamental domain of $\Pr^1\setminus\{a_1,\dots,a_{r_0}\}$
has area $1$.
When $D_g=r_0=3$, this area with respect to $\omega_{\Omega}$ 
is $G-1+3/2=1/2$ by (\ref{hyperbolicarea}), 
thus we use the metric $2\omega_{\Omega}$ in (\ref{Simizu2}), 
and we get 
\[
\begin{array}{ll}
T_g(r)\geq\displaystyle\int_0^r\dfrac{dt}{t}\int_{\C(t)}2g^*\omega_{\Omega}
=2\displaystyle\int_0^r\dfrac{dt}{t}\int_{\C(t)}\omega_h\\
=\log\dfrac1{1-r}\,.
\end{array}
\]
The last assertion in (i) follows from Proposision 5.3.
When $D_g=r_0=4$, we need no change of the metric, and get (ii).
\end{proof}
What we have to do is, however, the converse. 
That is, to estimate $T_g(r)$ first,
and then get a bound for $D_g$ or $\nu_{g}$. 
We leave this to \cite{KKM}. \\
\newline
{\bf Acknowledgement.} The authors thank the referee for a careful reading of the paper and adequate comments.

\end{document}